\documentclass[12pt,a4paper]{amsart}

\usepackage[dvips]{graphicx}
\newtheorem{theorem}{Theorem}[section]

\newtheorem{lemma}[theorem]{Lemma}

\theoremstyle{definition}

\theoremstyle{remark}
\newtheorem{remark}[theorem]{Remark}

\numberwithin{equation}{section}

\begin{document}
\title{Links with surgery yielding the 3-sphere}
\author[M. Teragaito]{Masakazu Teragaito}
\address{Department of Mathematics and Mathematics Education,
Hiroshima University, Kagamiyama 1-1-1, Higashi-hiroshima 739-8524, Japan}
\email{teragai@hiroshima-u.ac.jp}
%    General info

\begin{abstract}
For any $n\ge 2$,
we give infinitely many unsplittable links of $n$ components in the 3-sphere
which admit non-trivial surgery yielding the 3-sphere again and 
whose components are mutually distinct hyperbolic knots.
Berge and Kawauchi gave $2$-component hyperbolic links with those two properties.
We can also give infinitely many $2$-component hyperbolic links
with such properties.
\end{abstract}
\maketitle

\section{Introduction}

In \cite{GL}, Gordon and Luecke showed that non-trivial Dehn surgery on a non-trivial knot
in the $3$-sphere $S^3$ never yield $S^3$, and this led to the affirmative solution
of the knot complement conjecture.
On the other hand, there are abundant examples of links that admit
non-trivial surgeries yielding $S^3$.
For example, if a link has an unknotted componet, then $1/n$-surgery along the component
does not change the total space $S^3$ for any $n$.
Also, if there are two componets which cobound an annulus disjoint from
the other components of the link, then a suitable surgery along those two components,
corresponding to a Dehn twist along the annulus, gives $S^3$ again.
In \cite{B}, Berge gave a tunnel number one, hyperbolic link of two components
which admits three non-trivial surgeries yielding $S^3$ and
whose componets are distinct hyperbolic knots.
Kawauchi \cite{K} also announced the existence of such links shown
by his imitation theory.
The purpose of this paper is to give links with any number of components
which admit non-trivial surgery yielding $S^3$ and whose components are mutually
distinct hyperbolic knots by using Kirby calculus.
(See \cite{GS,R} for Kirby calculus.)
In particular, in case of two components, our links can be chosen to be tunnel number one and hyperbolic.
(Our examples are different from Berge's one.)

Let $L=K_1\cup \dots\cup K_n$ be a link of $n$ components.
For $\gamma=\gamma_1\cup\dots\cup \gamma_n$, where $\gamma_i$ is a slope on
$\partial N(K_i)$, we use $L(\gamma)$ (or $L(\gamma_1,\dots,\gamma_n)$) to denote the closed manifold obtained by
$\gamma$-surgery on $L$.
More precisely, $L(\gamma)$ is obtained by attaching $n$ solid tori $V_1,\dots,V_n$
to $S^3-{\mathrm{Int}}\,N(L)$ along their boundaries so that
each $\gamma_i$ bounds a meridian disk in $V_i$.
In this paper, we say that the surgery is \textit{non-trivial\/} if all $\gamma_i$ are 
non-meridional slopes.

\begin{theorem}\label{main}
For any $n\ge 2$, there exist infinitely many links $L=K_1\cup\dots\cup K_n$ of $n$ components
in $S^3$ satisfying the following:
\begin{itemize}
\item[(i)] $L$ admits a non-trivial surgery yielding $S^3$;
\item[(ii)] all $K_i$ are mutually distinct hyperbolic knots;
\item[(iii)] $L$ is unsplittable;
\item[(iv)] $L$ has tunnel number $n-1$.
\end{itemize}
In particular, when $n=2$, $L$ can be chosen to be hyperbolic.
\end{theorem}

%%%%%%%%%%%%%%%%%%%%%%%%%%%%%%%
\section{Proof of Theorem \ref{main}}

\begin{figure}[tb]
%\blankbox{1.0\columnwidth}{2in}
\includegraphics*[scale=0.3]{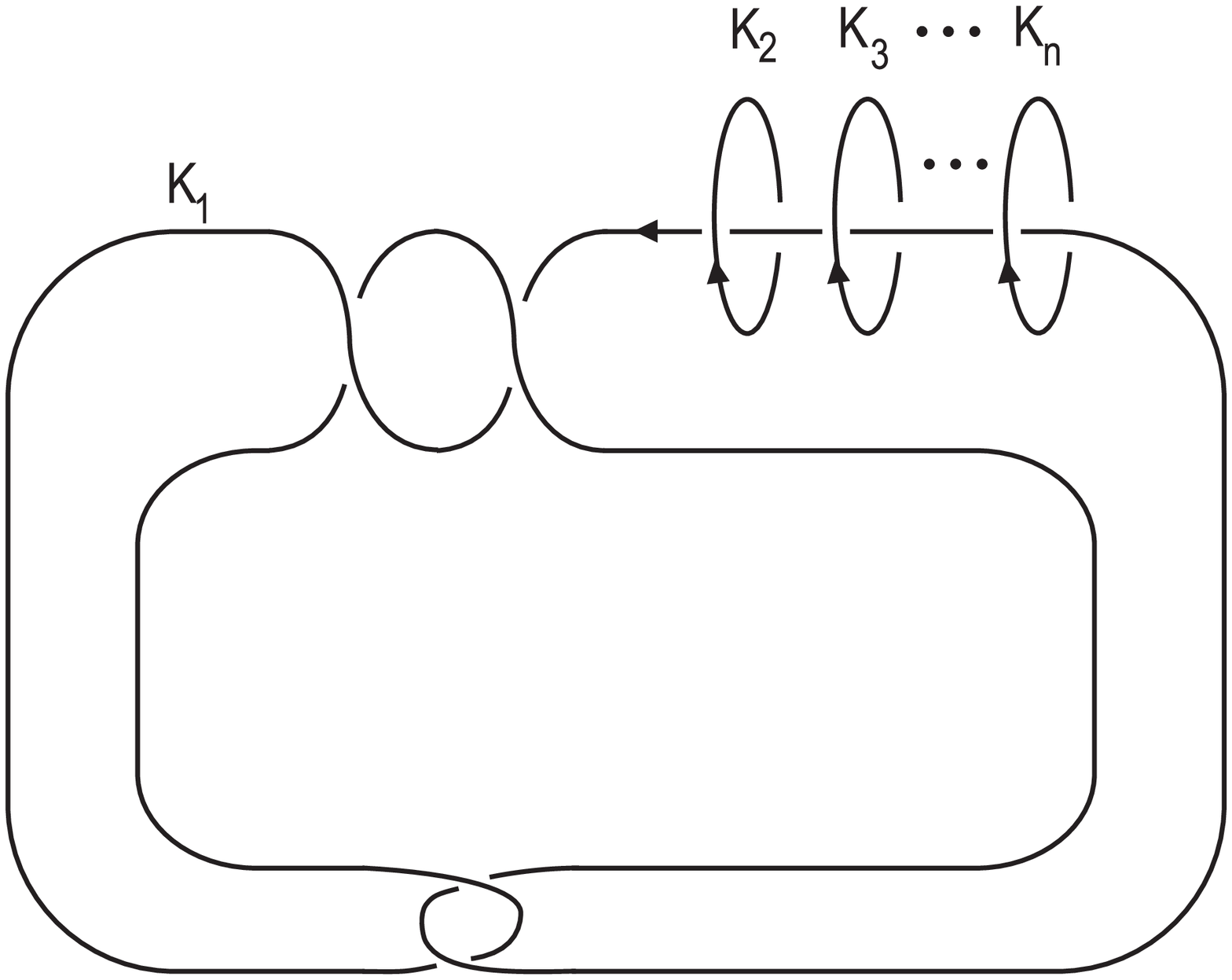}
\caption{}\label{start}
\end{figure}

Let $n\ge 2$ and $k\ge 0$ be integers.
Let $K_1$ be the figure-eight knot, and
let $L=K_1\cup K_2\cup\dots\cup K_n$, where $K_2,\dots,K_n$ are meridians of $K_1$.
Assign an orientation to $K_i$ as in Figure \ref{start}.

\begin{lemma}\label{base}
$L(n-2,0,1,1,\dots,1)=S^3$.
\end{lemma}

\begin{proof}
Perform $(-1)$-twists on $K_2,\dots,K_n$, and 
change $K_1$ to an unknot by using $K_2$.
Then the resulting framed link is a Hopf link
whose coefficients are $0$ and $2$.
It is well known that this presents $S^3$.
\end{proof}

\begin{figure}[tb]
%\blankbox{1.0\columnwidth}{2in}
\includegraphics*[scale=0.3]{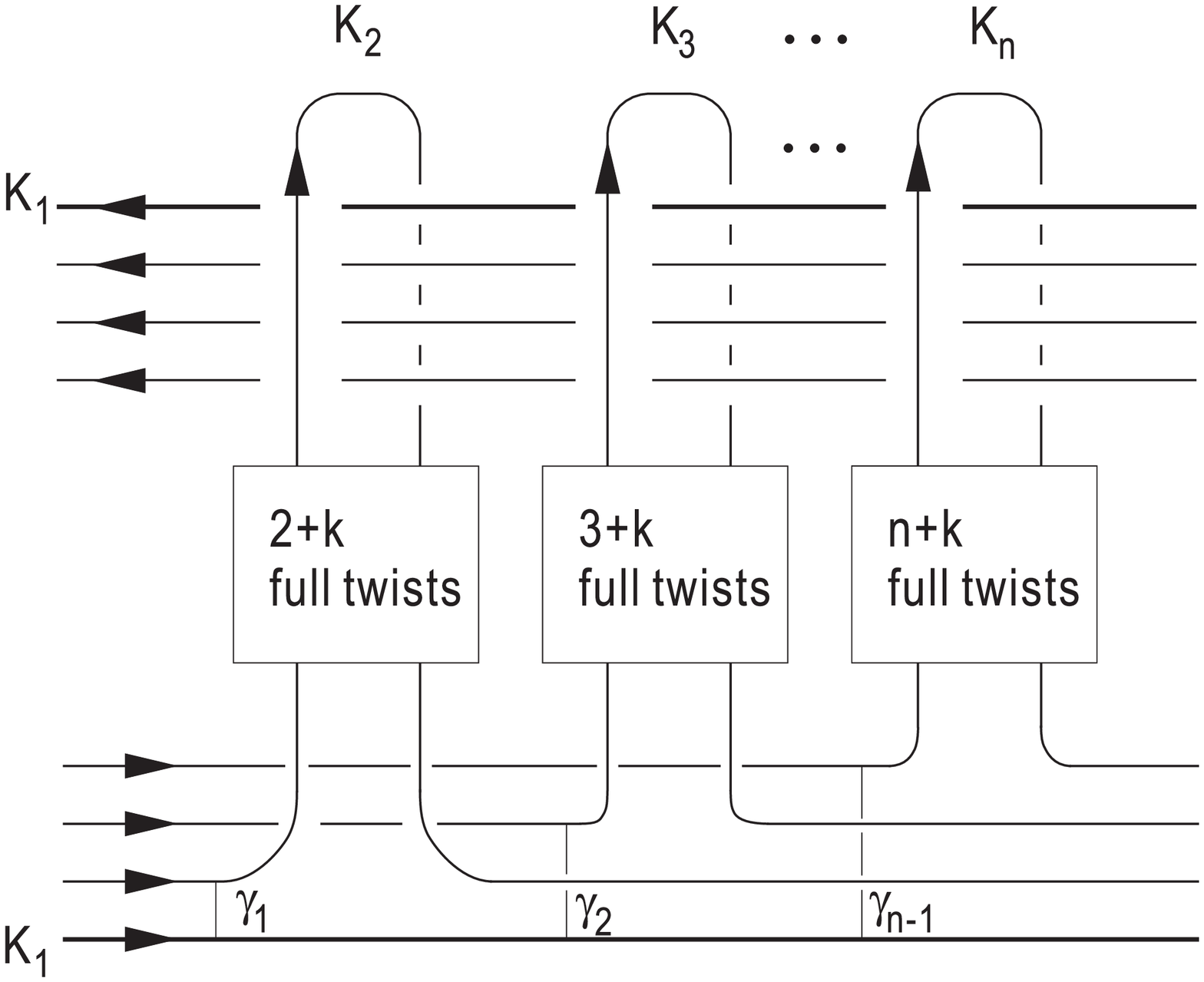}
\caption{}\label{slide}
\end{figure}

Now, slide $K_2,\dots,K_n$ over $K_1$.
More precisely, take a parallel copy $\ell_i$ of $K_1$ determining the framing on $K_1$ for $i=2,\dots,n$,
and form a band-sum of $K_i$ and $\ell_i$ as shown in Figure \ref{slide}, where
the band connecting $K_i$ and $\ell_i$ has $i+k$ full twists (right-handed).
Since $K_i$ and $K_1$ have linking number $1$, the resulting component, called $K_i$ again,
has coefficient $n$ if $i=2$, or $n+1$ if $i>2$.

This new link is also called $L$, and we now show that $L$ has the
desired properties.
By Lemma \ref{base}, $L(n-2,n,n+1,\dots,n+1)=S^3$.
It is straightforward to check that $K_i$ is the $2$-bridge knot $S(1+20(i+k),2-10(i+k))$.
Thus all $K_i$ are mutually distinct hyperbolic knots.
Since $K_i$ $(i\ge 2)$ and $K_1$ have linking number $n-1$, $L$ is unsplittable.
It is easy to see that
the arcs $\gamma_1,\dots,\gamma_{n-1}$ shown in Figure \ref{slide} form an unknotting tunnel system of $L$.
Finally, two different values of $k$ give distinct links.

When $n=2$, our link $L$ is tunnel number one.
Since each component of $L$ is hyperbolic $2$-bridge knot, $S^3-L$ is atoroidal by \cite{EU} and not Seifert fibered
by \cite{BM}.
Hence $L$ is hyperbolic by \cite{Th}.

This completes the proof of Theorem \ref{main}.

\begin{remark}
Even if $n>2$, our link $L$ seems to be hyperbolic, but
we could not prove it.  
\end{remark}

%%%%%%%%%%%%%%%%%%%%%%%%%%%%%%%%%%%%%%%%%%%%%%%%%%%%%%%%%%%%%%%%%%%

\bibliographystyle{amsplain}

\begin{thebibliography}{EU}

\bibitem{B} J. Berge,
\textit{Embedding the exteriors of one-tunnel knots and links in the 3-sphere},
unpublished manuscript. 

\bibitem{BM} G. Burde and K. Murasugi,
\textit{Links and Seifert fiber spaces},
Duke Math. J. \textbf{37} (1970), 89--93. 

\bibitem{EU} M. Eudave-Mu\~noz and Y. Uchida,
\textit{Non-simple links with tunnel number one},
Proc. Amer. Math. Soc. \textbf{124} (1996), 1567--1575. 

\bibitem{GS} R. Gompf and A. Stipsicz,
\textit{4-Manifolds and Kirby calculus},
Graduate Studies in Mathematics \textbf{20},
American Mathematical Society, Providence, RI (1999).

\bibitem {GL} C. McA. Gordon and J. Luecke,
\textit{Knots are determined by their complements},
J. Amer. Math. Soc. \textbf{2} (1989), 371--415.

\bibitem{K} A. Kawauchi,
private communication.

\bibitem{R} D. Rolfsen,
\textit{Knots and links},
Mathematics Lecture Series \textbf{7},
Publish or Perish, Inc., Berkeley, Calif. (1976). 

\bibitem{Th} W. Thurston,
\textit{Three dimensional manifolds, Kleinian groups and hyperbolic geometry},
Bull. Amer. Math. Soc. \textbf{6} (1982), 357--381.


\end{thebibliography}

\end{document}